\documentstyle[12pt,amssymbols,epsf]{article}

\def\presentation{
\voffset -.50in
\hoffset -.19in
\oddsidemargin 0in \evensidemargin 0in
\marginparwidth .75in \marginparsep 7pt \topmargin 0in
\headheight 12pt \headsep .25in
\footheight 18pt \footskip .35in
\textheight 9.5in \textwidth 6.5in
\columnsep 10pt \columnseprule 0pt }

\presentation

\newcommand{\beq}{\begin{equation}}
\newcommand{\eeq}{\end{equation}}
\newcommand{\bea}{\begin{eqnarray}}
\newcommand{\bean}{\begin{eqnarray*}}
\newcommand{\eea}{\end{eqnarray}}
\newcommand{\eean}{\end{eqnarray*}}
\newcommand{\demi}{\frac{1}{2}}

\newcommand{\Tr}{\mbox{Tr}}

\title{ A Quasi-Hopf algebra interpretation of quantum
3-j and 6-j symbols and difference equations.}
\author{O. Babelon $^{\dag} $ \and D. Bernard \thanks{Member of CNRS, Service
de
Physique Th\'eorique, CEN-Saclay, Orme des Merisiers,
F-91191, Gif-sur-Yvette, France and IHES 35 route de Chartres, 91440 Bures sur
Yvette, France.} \and E. Billey
\thanks{L.P.T.H.E. Universit\'e Paris VI (CNRS UA 280), Box 126, Tour 16,
$1^{{\rm er}}$ \'etage, 4 place Jussieu, 75252 Paris Cedex 05, France}}
\date{}

\begin{document}

\begin{titlepage}
\renewcommand{\thepage}{}
\maketitle
\vspace{2cm}
\begin{abstract}
We consider the universal solution of the Gervais-Neveu-Felder equation
in the ${\cal U}_q(sl_2)$ case. We show that it has a quasi-Hopf
algebra interpretation. We also recall its relation to quantum 3-j and
6-j symbols. Finally, we use this solution to build a q-deformation of
the trigonometric Lam\'e equation.
\end{abstract}

\vfill

PAR LPTHE 95-51, IHES/P/95/91

\end{titlepage}
\renewcommand{\thepage}{\arabic{page}}

\section{Introduction}
The Gervais-Neveu-Felder equation is a deformation of the standard
Yang-Baxter equation. In the $sl_2$ case, it reads
 \begin{eqnarray}
 R_{12}(x) R_{13}(x q^{H_2}) R_{23}(x) =
R_{23}(x q^{H_1}) R_{13}(x) R_{12}(x q^{H_3})
\label{GNF}
\end{eqnarray}
Here, $H$ denotes a Cartan generator in $sl_2$
(or rather ${\cal U}_q(sl_2)$) and  $x$ is a
parameter not to be confused with the spectral parameter
(absent in the $sl_2$ case).

This equation appeared independently in several contexts. It was first
discovered by J.L. Gervais and A. Neveu in their studies on Liouville
theory \cite{GeNe}. It
was rediscovered by G. Felder in his approach to the quantization of the
Knizhnik-Zamolodchikov-Bernard equation \cite{Fel1}. Finally, it was shown to
play an
important role in the quantization of the Calogero-Moser models in the
$R$-matrix formalism \cite{ABB95}. For all these reasons, we believe that this
equation deserves much
attention.

In this note, we analyse the universal solution $R(x) \in {\cal
U}_q(sl_2) \otimes {\cal U}_q(sl_2)$ of eq.(\ref{GNF}) obtained in
\cite{Ba1}. We show that it has a nice quasi-Hopf algebra
interpretation. For completeness, we recall the connection of this
solution with $q$-analogs of 3-j and 6-j symbols. Finally, we explain
how it can be used to construct a $q$-difference analog of the
trigonometric Lam\'e equation (Calogero model for 2 particules).

\section{A summary of universal formulae}

In this section, we recall the universal formulae obtained in \cite{Ba1} for
the matrix $R_{12}(x) \in {\cal U}_q(sl_2) \otimes {\cal U}_q(sl_2)$.
We denote by $H,E_\pm$ the generators of the quantum group ${\cal U}_q(sl_2)$
\begin{eqnarray}
 \relax [ H,E_\pm] = \pm 2 E_\pm,  ~~~~
 \relax [ E_+ , E_- ] = { q^H - q^{-H} \over q - q^{-1} } \nonumber
 \end{eqnarray}
The coproduct is defined as
\begin{eqnarray}
\Delta (H) = H\otimes id + id \otimes H,~~~~
\Delta (E_\pm) = E_\pm \otimes q^{{1\over 2} H } +
   q^{- {1\over 2} H } \otimes  E_\pm \nonumber
   \end{eqnarray}
We have $R^D_{12} \Delta(a) = \Delta'(a) R^D_{12}$
for any $a\in U_q(sl_2)$  where $\Delta'$ is the opposite
comultiplication and
$R^D_{12}$ Drinfeld's universal $R$-matrix~:
\begin{eqnarray}
R^D_{12} = q^{{1\over 2} H \otimes H} \sum_{i=0}^\infty (q -q^{-1})^i \;
{q^{-{i(i+1)\over 2}} \over [i]! }\; q^{{i\over 2}H} E_+^i \otimes
q^{-{i\over 2}H} E_-^i
\nonumber
\end{eqnarray}
As usual, q-numbers are defined as $[i] = (q^i -q^{-i})/(q-q^{-1})$.
Let us now define
\begin{eqnarray}
 R_{12}(x) = F_{21}^{-1}(x)\; R^D_{12}\; F_{12}(x)
  \label{rshift}
   \end{eqnarray}
with
\begin{eqnarray}
F_{12}(x)& =& \sum_{k=0}^\infty (q-q^{-1})^k {(-1)^k \over [k]! }
 {x^k \over \prod_{\nu =k}^{2k-1} (x q^\nu q^{H_2} - x^{-1}q^{-\nu} q^{-H_2} )}
 q^{{k\over 2}(H_1 +H_2)} E_+^k \otimes E_-^k
 \label{fdex} \\
F_{12}^{-1}(x) &=& \sum_{k=0}^\infty (q-q^{-1})^k {1 \over [k]! }
 {x^k \over \prod_{\nu =1}^{k} (x q^\nu q^{H_2} - x^{-1}q^{-\nu}
 q^{-H_2}
 )}q^{{k\over 2}(H_1 +H_2)} E_+^k \otimes E_-^k \nonumber
 \end{eqnarray}
It follows from the construction of \cite{Ba1} that $R_{12}(x)$ is a solution
of eq.(\ref{GNF}).

One can check that $F_{12}(x)$ satisfies the
following ``shifted cocycle'' condition
\begin{eqnarray}
\left[ (id \otimes \Delta )F \right] \cdot
\left[ id \otimes F \right] =
\left[  (\Delta \otimes id )F \right] \cdot \left[
F(xq^{H_3})\otimes id \right]
\label{cocycle}
\end{eqnarray}
This relation is proved using standard q-binomial identities.
It turns out that $F_{12}(x)$ is actually a ``shifted coboundary''
\begin{eqnarray}
F_{12}(x) = \Delta M(x)~\Big[ id \otimes  M(x)\Big]^{-1}
\Big[ M(xq^{H_2}) \otimes id \Big]^{-1}
\label{cobord}
\end{eqnarray}
where the formula for the ``boundary'' reads
\begin{eqnarray}
M(x) = \sum_{n,m=0}^\infty { (-1)^m x^m q^{{1\over 2}n(n-1) + m(n-m)}
\over [n]! [m]! \prod_{\nu =1}^n (x q^\nu - x^{-1} q^{-\nu})}
E_+^n E_-^m q^{{1\over 2}(n+m)H}
\nonumber
\end{eqnarray}
Equation (\ref{cobord}) implies eq. (\ref{cocycle}).

\section{Quasi-Hopf algebra interpretation}

The previous construction possesses a natural
quasi-Hopf interpretation.
Indeed, since $R(x)$
is defined in eq.(\ref{rshift}) by a twisting procedure
in the sense of Drinfeld \cite{drquasi} it is canonically
associated to a quasi-Hopf structure on $U_q(sl_2)$.
We shall denote it as $U_{q;x}(sl_2)$.
This quasi-Hopf algebra possesses very specific properties
due to the ``shifted cocycle" relation
(\ref{cocycle}) satisfied by $F(x)$. Besides Drinfeld's construction, this
gives
another example of a quasi-Hopf algebra structure over ${\cal U}_q(sl_2)$.

Let us recall following ref.\cite{drquasi} that
a quasi-Hopf algebra is specified by a
quadruplet $(A,\Delta, R, \Phi)$ where $A$ is an
assocative algebra, $\Delta$ is  a (non-coassociative)
comultiplication in $A$, $R\in A\otimes A$ and
$\Phi \in A\otimes A \otimes A$ are such that~:
\begin{eqnarray}
R \Delta(a) &=& \Delta'(a) R \label{quasiR}\\
(id\otimes \Delta) \Delta(a) ~ \Phi
&=& \Phi~ (\Delta\otimes id) \Delta(a) \label{quasiD}
\end{eqnarray}
for all $a\in A$. There also are extra compatibility
relations between $\Delta,~R$ and $\Phi$ which we
shall mention when needed. We will consider quasitriangular
quasi-Hopf algebra, i.e., $R$ is assumed
to verify the conditions
\begin{eqnarray}
(\Delta \otimes id )R &=& \Phi_{321} R_{13} \Phi^{-1}_{132} R_{23} \Phi_{123}
\nonumber
\\
(id \otimes \Delta )R &=& \Phi^{-1}_{231} R_{13} \Phi_{213} R_{12}
\Phi^{-1}_{123}
\nonumber
\end{eqnarray}

\def\[{\big[}
\def\]{\big]}
There exists a twisting procedure to
construct quasi-Hopf algebras. Namely,
if $(A,\Delta,R,\Phi)$ is a quasitriangular
quasi-Hopf algebra then
a new quasitriangular quasi-Hopf algebra
$(A,\widetilde \Delta, \widetilde R, \widetilde \Phi)$  is
defined by
$\widetilde \Delta(a)= F^{-1}_{12}\Delta(a) F_{12}$,
and
\begin{eqnarray}
\widetilde \Phi&=& F^{-1}_{23} \[{(id\otimes \Delta)(F^{-1})}\]
 ~\Phi ~ \[{ (\Delta\otimes id)(F)}\] F_{12} \label{phitwist}\\
\widetilde R &=& F^{-1}_{21} R F_{12} \label{rtwist}
\end{eqnarray}
with $F_{12}\in A\otimes A$.

In our case, we are twisting $U_q(sl_2)\equiv (U_q(sl_2), \Delta, R^D,id) $ by
$F(x)$.
So we have $\widetilde R = F^{-1}_{21}(x) R^D_{12} F_{12}(x)= R(x)$
as defined in eq.(\ref{rshift}). We denote
$\widetilde \Delta$  by $\Delta_x$ with~:
\begin{eqnarray}
\Delta_x (a) = F_{12}^{-1}(x) \Delta (a) F_{12}(x),
\quad \forall a\in U_q(sl_2)
\label{deltax}
\end{eqnarray}
It is a simple check to verify that the ``shifted cocycle"
condition (\ref{cocycle}) implies that~:
\begin{eqnarray}
(id \otimes \Delta_x) \Delta_x (a) =
(\Delta_{xq^{H_3}}\otimes id) \Delta_x (a)
\label{cox}
\end{eqnarray}
In other words, the shift breaks the co-associativity.
We denote $\widetilde \Phi$
by $\Phi(x)$. It possesses a simple expression
in terms of $F(x)$~:
\begin{eqnarray}
\Phi(x)&=& F^{-1}_{23}(x) \[{(id\otimes \Delta)(F^{-1}(x))}\]
\[{(\Delta\otimes id)(F(x))}\] F_{12}(x) \nonumber\\
 &=& F^{-1}_{12}(xq^{H_3})~ F_{12}(x) \label{phix}
 \end{eqnarray}
where in the last equality we again used the ``shifted cocycle"
relation (\ref{cocycle}).

We can now write all the quasi-Hopf relations in
$U_{q;x}(sl_2)$ in terms of $R(x)$ or $F(x)$.
For example, the general relation (\ref{quasiD})
reduces to eq.(\ref{cox}). Also, thanks
to the following property,
\begin{eqnarray}
R_{12}(xq^{H_3})= \Phi_{213}(x) R_{12}(x) \Phi^{-1}_{123}(x)
\nonumber
\end{eqnarray}
the quasi- Yang-Baxter equation,
\begin{eqnarray}
 \Phi_{321}^{-1}(x) R_{12}(x) \Phi_{312}(x) R_{13}(x) \Phi^{-1}_{132} R_{23}(x)
 =  R_{23}(x) \Phi^{-1}_{231}(x) R_{13}(x)
 \Phi_{213}(x) R_{12}(x) \Phi^{-1}_{123}(x) \nonumber
 \end{eqnarray}
valid in any quasitriangular quasi-Hopf algebra reduces to the equation
(\ref{GNF}).

Similarly, the quasitriangular property of $U_{q;x}(sl_2)$
implies that~
\begin{eqnarray}
(\Delta_x\otimes id) R(x) &=& R_{13}(xq^{H_2}) R_{23}(x)
F^{-1}_{12}(xq^{H_3}) F_{12}(x) \nonumber\\
(id \otimes \Delta_x) R(x) &=&
F^{-1}_{23}(x) F_{23}(xq^{H_1})
R_{13}(x) R_{12}(xq^{H_3}) \nonumber
\end{eqnarray}

Notice that for $x=0$, $F_{12}(x)\vert_{x=0}=1$ and therefore
$R(x)\vert_{x=0}= R^D$.
In the limit $x=\infty$,
$F^{-1}_{12}(x)\vert_{x=\infty} =  q^{-H\otimes H/2} R^D_{12}$. Thus
$R_{12}(x)\vert_{x=\infty} = q^{-H\otimes H/2} R^D_{21} q^{H\otimes H/2}$,
and $\Delta_{x=\infty}(a) = q^{-H\otimes H/2} \Delta'(a)
q^{H\otimes H/2}$ for all $a\in U_q(sl_2)$.

\section{Relation to 3-j and 6-j symbols}

We now give a list of formulae expressing the matrix elements of the various
objects we have considered so far in terms of standard q-analogs of the
3-j and 6-j symbols. Let $\rho^{(j)}$ denote the spin j representation of
${\cal U}_q(sl_2)$. Then
\begin{eqnarray}
\rho^{(j)}(H) \vert j,m \rangle &=& 2 m ~ \vert j,m \rangle \nonumber \\
\rho^{(j)}(E_\pm) \vert j,m \rangle &=& \sqrt{[j \mp m] [j \pm m +1] }
{}~\vert j,m \pm 1 \rangle \nonumber
\end{eqnarray}
The first step is to find the matrix elements of the matrix $M(x)$ in the
spin-j representation. We get
\begin{eqnarray}
\left[ M^{(j_1)}(x) \right]_{\sigma_1\, m_1}&=& (-1)^{\sigma_1+m_1}
{ \sqrt{ [j_1 +\sigma_1 ]![j_1 -\sigma_1 ]![j_1+m_1 ]! [j_1-m_1 ]!}\over
 \prod_{r=1}^{j_1+\sigma_1} (1-x^2 q^{2r} ) }\cdot \label{Mj1} \\
 &&\cdot q^{\sigma_1(\sigma_1-m_1)} x^{\sigma_1-m_1} \sum_p {q^{2p\,
 \sigma_1}  x^{2p}
 \over [p]!~ [\sigma_1-m_1 +p]! [j_1-\sigma_1-p]!~ [j_1+m_1 -p]! }
 \nonumber
\end{eqnarray}

This formula agrees (up to normalizations) with the one found in
\cite{Ger90}.

This matrix $M(x)$ is known to perform the vertex-IRF transformation in
conformal field theory \cite{Pas,Ba2,MooRes,FFK}.
\begin{eqnarray}
\xi_{m_1}^{(j_1)} (z) = \sum_{\sigma_1} \psi_{\sigma_1}^{(j_1)}(z)
 M^{(j_1)}_{\sigma_1 m_1}(x)
\nonumber
\end{eqnarray}
where the $\psi$'s are IRF type operators and the $\xi$'s are vertex
type operators. The braiding relations of the $\psi$'s are described by the
matrix $R(x)$, while those of the $\xi$'s are described by $R^D$.
 Thus, we expect the elements
$M^{(j_1)}_{\sigma_1m_1}(x)$ to be related to quantum 3-j symbols. The
precise connexion was found in \cite{CGR94}.  We have
\begin{eqnarray}
\left[ M^{(j_1)}(x) \right]_{\sigma_1\, m_1}&=& { {\cal N}^{(j_1)}_\psi
(x,\sigma_1)\over {\cal N}^{(j_1)}_\xi
(m_1) } \lim_{m \to \infty}  \left[ \matrix{j_1 & j(x) & j(x)+\sigma_1 \cr
m_1 & m & m+ m_1}
\right]_q
\label{M=3j}
\end{eqnarray}
where we have defined $j(x)$ through the relation
\begin{eqnarray}
x = q^{2j(x)+1}
\label{jx}
\end{eqnarray}
Eq.(\ref{M=3j}) has to be understood as an analytic continuation in $j(x)$ of
3-j
symbols \cite{GeRo}.
We give a sketch of the proof in the Appendix.  The factors
${\cal N}^{(j_1)}_\xi$ and
${\cal N}^{(j_1)}_\psi$ can be reabsorbed into the normalizations of the fields
$\psi$ and $\xi$
respectively. Their expression is also given in the Appendix. We represent
eq.(\ref{M=3j}) by a diagram

\begin{eqnarray}
 [M^{(j_1)}(x)]_{\sigma_1 m_1}=\raisebox{0 cm}{\mbox{\epsffile{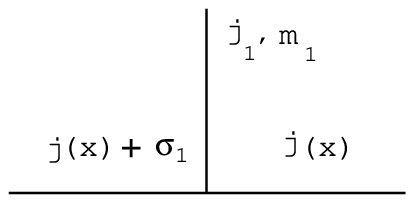} }}
\end{eqnarray}

{}From eq.(\ref{M=3j}), it is now possible to build the complete dictionary
between the matrix elements of $F_{12}(x)$ and $R_{12}(x)$ and standard
3-j and 6-j symbols.

We start with
\begin{eqnarray}
\langle j_1, \sigma_1\vert \langle j_2, \sigma_2\vert
 M_2^{(j_2)}(xq^{H_1}) M_1^{(j_1)}(x) \vert j_1,m_1 \rangle \vert
 j_2,m_2 \rangle  &=&  M_{\sigma_2,m_2}^{(j_2)}(xq^{2\sigma_1})
  M_{\sigma_1,m_1}^{(j_1)}(x)= \nonumber \\
  & & \nonumber \\
  && \hskip -10cm
 { {\cal N}^{(j_1)}_\psi (x,\sigma_1) {\cal N}^{(j_2)}_\psi (xq^{2\sigma_1},
 \sigma_2)
  \over {\cal N}^{(j_1)}_\xi(m_1) {\cal N}^{(j_2)}_\xi(m_2) }
  \lim_{m,m' \to \infty}  \left[ \matrix{j_2 & j(x)+\sigma_1
  & j(x)+ \sigma_1 + \sigma_2 \cr m_2 &
   m & m +m_2 }
  \right]
  \left[ \matrix{j_1 & j(x) & j(x)+\sigma_1 \cr
  m_1 &  m' & m'+m_1 }
  \right]
  \nonumber
  \end{eqnarray}
Notice that we have used the fact that
\begin{eqnarray}
xq^{2\sigma_1} = q^{2(j(x)+\sigma_1) +1}
\nonumber
\end{eqnarray}
Hence the shift in the Gervais-Neveu-Felder equation $x \to xq^H$ precisely
corresponds to the shift of spins $j(x) \to j(x) +\sigma $.
Thus we have the diagramatic correspondance
\begin{eqnarray}
M_2^{(j_2)}(xq^{H_1}) M_1^{(j_1)}(x) = \raisebox{-.0
cm}{\mbox{\epsffile{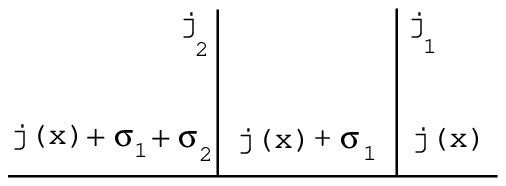} }}
\end{eqnarray}

The matrix elements of $R_{12}(x)$ are computed from the formula
\begin{eqnarray}
R_{12}(x) M_1(xq^{H_2})M_2(x) =  M_2(xq^{H_1}) M_1(x) R^D_{12}
\nonumber
\end{eqnarray}
or graphically
\begin{eqnarray}
\sum_{\sigma_1 \sigma_2} R(x)^{j_1j_2}_{\sigma'_1 \sigma'_2, \sigma_1 \sigma_2}
\raisebox{-1.2 cm}{\mbox{\epsffile{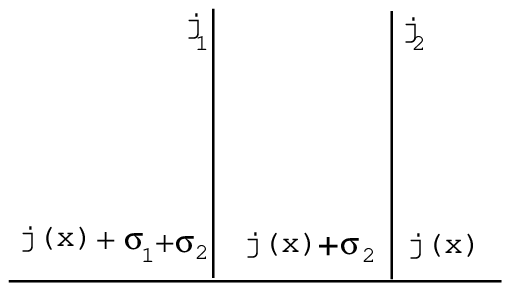} }} =  \raisebox{-1.2
cm}{\mbox{\epsffile{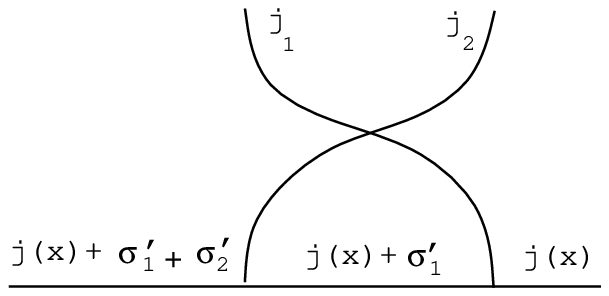} }}
\end{eqnarray}

This is equivalent to the braiding relation and relates the matrix elements
of $R(x)$ to 6-j symbols:
\begin{eqnarray}
\langle j_1, \sigma_1'\vert \langle j_2, \sigma_2' \vert
 R_{12}(x) \vert j_1,\sigma_1 \rangle \vert
  j_2,\sigma_2 \rangle  &=&  (-1)^{\sigma_1' -\sigma_1}
  q^{C(j(x))+C(j(x)+\sigma_1 +\sigma_2) -
  C(j(x)+\sigma_1')-C(j(x)+\sigma_2) }
  \nonumber \\
&&\hskip -5cm { {\cal N}^{(j_1)}_\psi (x,\sigma_1') {\cal N}^{(j_2)}_\psi
(xq^{2\sigma_1'}
,\sigma_2')
\over {\cal N}^{(j_1)}_\psi (xq^{2\sigma_2},\sigma_1){\cal N}^{(j_2)}_\psi
(x,\sigma_2)}
\left\{ \matrix{j_2 & j(x)+\sigma_1 +\sigma_2 & j(x)+\sigma_1' \cr
  j_1 & j(x) & j(x) +\sigma_2 } \right\}_q
\nonumber
\end{eqnarray}
where $C(j) = j(j+1)$ and the last symbol represents a 6-j coefficient
(see eq. 5.11 in \cite{KiResh}).

Finally, we give the formula for the matrix elements of $F_{12}(x)$
in terms of 3-j and 6-j symbols. We start from the formula
\begin{eqnarray}
F_{12}(x) M_1(xq^{H_2}) M_2(x) = \Delta M(x)
\nonumber
\end{eqnarray}
{}From the definition of the coproduct, we have
\begin{eqnarray}
\relax[ \Delta M^{j_1 j_2}(x) ]_{\sigma_1 \sigma_2 , m_1 m_2} =
\sum_{j_{12}}  \left[\matrix{ j_1 & j_2 & j_{12} \cr \sigma_1 & \sigma_2 &
\sigma_1+\sigma_2 } \right]_q \Big[ M^{(j_{12})}(x)\Big]_{\sigma_{1}+\sigma_2,
m_{1}+m_2}
\left[\matrix{ j_1 & j_2 & j_{12} \cr m_1 & m_2 & m_1+m_2 } \right]_q
\nonumber
\end{eqnarray}
Using the interpretation of $M$ as a 3-j symbol and the defining relation of
6-j symbols,
we can write
\begin{eqnarray}
 \Big[ M^{(j_{12})}(x)\Big]_{\sigma_{1}+\sigma_2, m_{1}+m_2}
 \left[\matrix{ j_1 & j_2 & j_{12} \cr m_1 & m_2 & m_{1}+m_2 } \right]_q&=&
\nonumber \\
&&\hskip -9cm \sum_{\sigma_1' \sigma_2'}{{\cal N}^{(j_{12})}_\psi (x,\sigma_{1}
+\sigma_2) \over
 {\cal N}^{(j_1)}_\psi (xq^{2\sigma_2'},\sigma_1') {\cal N}^{(j_2)}_\psi
 (x,\sigma_2') }
 \left\{ \matrix{ j_1 & j_2 & j_{12} \cr
   j(x) & j(x) +\sigma_{1}+\sigma_2 & j(x) + \sigma_2' } \right\}_q
   M^{(j_1)}_{\sigma_1',m_1}(xq^{2\sigma_2'})
   M^{(j_2)}_{\sigma_2',m_2}(x)
   \nonumber
\end{eqnarray}
Hence
\begin{eqnarray}
\relax[F^{j_1 j_2}(x) ]_{\sigma_1 \sigma_2 , \sigma_1' \sigma_2'}
&=& \nonumber \\
&& \hskip -5cm\sum_{j_{12}} {{\cal N}^{(j_{12})}_\psi (x,\sigma_{1}+\sigma_2)
\over
 {\cal N}^{(j_1)}_\psi (xq^{2\sigma_2'},\sigma_1') {\cal N}^{(j_2)}_\psi
  (x,\sigma_2') }
  \left[\matrix{ j_1 & j_2 & j_{12} \cr \sigma_1 & \sigma_2 &
\sigma_1+\sigma_2 } \right]_q
 \left\{ \matrix{ j_1 & j_2 & j_{12} \cr
    j(x) & j(x) +\sigma_1'+\sigma_2' & j(x) + \sigma_2' } \right\}_q
\nonumber
\end{eqnarray}

\section{Application to the trigonometric q-deformed Lam\'e equation}

In \cite{ABB95} it was shown how solutions of eq.(\ref{GNF}) could be used to
construct a set of commuting operators
corresponding to $q$-deformations of the quantum Calogero-Moser
Hamiltonians. In the ${\cal U}_q(sl_2)$ case, there is only one such
operator once we separate the center of mass motion.

According to the general prescription \cite{ABB95}, we start from a Lax matrix
satisfying
\beq
\label{Lalgebra}
R_{12}(x q^{-\demi H_3}) L_{13}(x) L_{23}(x) = L_{23}(x) L_{13}(x) R_{12}
(x q^{\demi H_3}),
\eeq
with a subscript 3 denoting the quantum space. The first
Hamiltonian is the restriction of $\Tr_1(L_{13}(x))$ to the subspace
of zero-weight vectors.

In the representation $\rho=\rho^{(1/2)} \otimes \rho^{(1/2)}$
of ${\cal U}_q(sl_2) \otimes {\cal U}_q(sl_2)$, the following extra
property is true:
$$ \rho \left ( [(H_1+H_2)\partial_x,R_{12}(x)] \right ) =0.$$
This condition allows to recast eq.(\ref{GNF}) in the form
(\ref{Lalgebra}) with a Lax operator $L(x)$ obtained by dressing
$R(x)$ with suitable shift operators:
$$ L_{13}(x)=q^{-(H_1+\demi H_3)p} \ R_{13}(x) \ q^{\demi H_3 p},
\quad {\rm with} \quad p=x {\partial\over \partial x}.$$

In the representation $\rho_j=\rho^{(1/2)} \otimes \rho^{(j)}$ of
${\cal U}_q(sl_2) \otimes {\cal U}_q(sl_2)$,
\beq
\label{Lj}
\rho_j(L_{13}(x))=\left (
\begin{array}{c c}
q^{-p} q^{\demi H} &  - q^{- \demi} x^{-1} f(x q^{\demi H})
q^{-\demi H} E_- \\
 & \\
q^{- \demi} x f(x q^{-\demi H+1}) q^{\demi H} E_+ &
q^p q^{-\demi H} \left [ 1 - f(x q^{-\demi H}) f(x q^{\demi H-1})
E_+ E_- \right ]
\end{array}
\right )
\eeq
with $f(x)=(q-q^{-1}) / (x-x^{-1}).$

Taking the trace on the first space we get
$$ \Tr_1(L_{13}(x)) = q^{-p} q^{\demi H} + q^p q^{-\demi H} \left [
1 - f(x q^{-\demi H}) f(x q^{\demi H-1}) E_+ E_- \right ].$$
We still have to restrict this operator to the space of zero-weight
vectors. In the spin $j$, representation, when $j$ is
integer, this subspace is one-dimensional
and the resulting Hamiltonian is scalar. Using $E_+ E_- \vert j,0
\rangle = [j] [j+1] \vert j,0 \rangle$, we get
$$ H_j = q^{-p} + q^p \left ( 1 - \frac{(q-q^{-1})^2 [j] [j+1]}
{(x-x^{-1})(q^{-1} x- q x^{-1})} \right ).$$
At the first non-trivial
order of $H_j$ in the limit $q \rightarrow 1$, we recover
the Calogero-Moser Hamiltonian
$ - \partial_z^2 + j (j+1) / \sinh^2(z)$, with  $x=\exp(z).$ Notice that
the coupling constant is related to the spin of the representation.

Alternatively, introducing the function
$$ c_j(x)=\frac{(q^j x - q^{-j} x^{-1})
(q^{-j-1} x - q^{j+1} x^{-1})}{(x-x^{-1})(q^{-1} x- q x^{-1})},$$
the Hamiltonian $H_j$ is given by
$$ H_j = q^{-p} + q^p c_j(x).$$
The eigenfunctions $\Psi$ of $H_j$ are the solutions
of the following trigonometric $q$-deformed Lam\'e equation:
$$ \Psi(q^{-1}x) + c_j(qx) \Psi(q x) = E \ \Psi(x).$$
An elliptic version of this equation already appeared in a different
context in \cite{KriZa}.

Integrability of the system manifests itself in the  fact that we
can easily solve this equation by using, for instance, the following
recursive procedure. For $j=0$ the Hamiltonian $H_0 = q^p + q^{-p}$ is free;
its eigenfunctions
are plane waves $\Psi_0(x,k) = x^k$ with the corresponding energy
$E(k)=q^k+q^{-k}.$ Let us now introduce the following ``shift'' operator
$$\quad D_j = q^{-p} - q^p \frac{(q^{-j} x - q^j x^{-1})
(q^{-j-1} x - q^{j+1} x^{-1})}{(x - x^{-1}) (q^{-1} x - q x^{-1})}$$
which satisfies
$$ H_j D_j = D_j H_{j-1},$$
The eigenfunction $\Psi_j(x,k)$ of
$H_j$ with energy $E(k)=q^k + q^{-k}$, are obtained by the successive action of
the
``shift'' operator:
$$ \Psi_j(x,k) = D_j \ \Psi_{j-1}(x,k).$$

Since the energy is even in $k$, we can start the recursion with $\Psi_0(x,k) =
x^k-x^{-k}$.
Then we get
\begin{eqnarray}
 \Psi_j(x,k)&=& \sum_{n=0}^j (-1)^n \left[ \matrix{j \cr n}\right]_q
 {\prod_{r=1}^n (q^{r-j-1}x - q^{-r+j+1}x^{-1})\over
  \prod_{r=1}^n (q^{r}x - q^{-r}x^{-1})}(q^{k(2n-j)}x^k - q^{-k(2n-j)}x^{-k})
\nonumber
\end{eqnarray}

This wave function has the interesting properties that the residues at the
poles
$x=\pm q^{-r}$ for $1 \leq r \leq j$ all vanish, and moreover, one has
$\Psi_j(x,k)=0$ for
$k= -j, -j+1,\cdots, j$. This is an analogue of the generalized
exclusion principle present in the Calogero-Sutherland model
\cite{BGHP}.

\section{Appendix}

We give an idea of the proof of eq.(\ref{M=3j}). We adopt here a naive
point of vue. We refer to \cite{CGR94} for a more detailed discussion.
We start from an avatar of van der Waerden formula for 3-j symbols (combine
eq. 3.5 and eq. 3.10 in ref.\cite{KiResh}):
\begin{eqnarray}
\relax \left[ \matrix{j_1 & j_2 & j_3 \cr m_1 & m_2 & m_3 } \right]_q &=&
\delta_{m_1+m_2 ,m_3}~ \Delta (j_1,j_2,j_3)~ q^{-{1\over 2} (j_1+j_2
-j_3)(j_1+j_2+j_3+1)
+j_1m_2 - j_2 m_1} \sqrt{[2j_3+1]}\cdot \nonumber \\
&&\hskip 1cm \cdot  \sqrt{ [j_1+m_1]!  [j_1-m_1]!  [j_2+m_2]!
 [j_2-m_2]!  [j_3+m_3]!  [j_3-m_3]! }\cdot \nonumber \\
 &&\hskip -3cm \cdot \sum_p { (-1)^p q^{p(j_1+j_2 +j_3 +1)} \over
 [p]! [j_1+j_2-j_3-p]! [j_2-m_2-p]! [j_1+m_1 -p]! [j_3-j_1+m_2+p]!
[j_3-j_2-m_1+p]! }
 \nonumber
 \end{eqnarray}
where
\begin{eqnarray}
\Delta (j_1,j_2,j_3) &=&(-1)^{j_1+j_2-j_3} \sqrt{ { [-j_1+j_2+j_3]!
[j_1-j_2+j_3]! [j_1+j_2 -j_3]!
\over [j_1+j_2+j_3 +1]!} }
\nonumber
\end{eqnarray}

We take a limit $m_2 \to \infty$ such that
\begin{eqnarray}
\lim_{m_2 \to \infty} q^{m_2} = 0, \quad \lim_{m_2 \to \infty} q^{-m_2} =
\infty
\nonumber
\end{eqnarray}
Then, one has
\begin{eqnarray}
 { [\alpha \pm m_2 ]! \over [\beta \pm m_2 ]! } \sim (\mp)^{\alpha -\beta}
{q^{\mp {1\over 2} (\alpha -\beta)(\alpha +\beta +1) } \over (q-q^{-1})^{\alpha
-
\beta} } q^{-(\alpha -\beta)m_2}
\nonumber
\end{eqnarray}
To perform the limit, we write the terms containing $m_2$ in the following form

\begin{eqnarray}
 \sqrt{{ [j_2+m_2]! \over  [j_3-j_1 +m_2]!}\cdot
{[j_3+m_1+m_2]!\over  [j_3-j_1 +m_2]! }
 \cdot {[j_2-m_2]!\over [j_2-m_2]!}\cdot { [j_3-m_1 -m_2]!  \over [j_2 -m_2]! }
}
 &&\sim \nonumber \\
 && \hskip -13cm (-1)^{j_1+{1\over 2}(j_2-j_3+m_1)}~
 { q^{ {1\over 2}(-2j_3m_1-m_1-2j_1(j_3+1)
 +j_1(j_1+1)+j_3(j_3+1) -j_2(j_2+1))} \over (q-q^{-1})^{j_1} }~q^{-j_1m_2}
\nonumber
\end{eqnarray}
and
\begin{eqnarray}
\lim_{m_2 \to \infty} { [j_2 -m_2]!~ [j_3-j_1 +m_2]! \over  [j_2-m_2-p ]!
[j_3-j_1+m_2 +p]! }
= (-1)^p q^{p(j_2+j_3-j_1+1)}
\nonumber
\end{eqnarray}
This decomposition is to ensure that we get the above important sign $(-1)^p$
correctly. Hence
\begin{eqnarray}
\lim_{m_2 \to \infty}  \left[ \matrix{j_1 & j_2 & j_3 \cr m_1 & m_2 & m_1+m_2 }
\right]_q &=&
\Delta (j_1,j_2,j_3)~{
\sqrt{[2j_3+1]}\sqrt{ [j_1+m_1]! [j_1-m_1]!  }\over (q-q^{-1})^{j_1}}\cdot
\nonumber \\
 &&\cdot (-1)^{j_1+{1\over 2}(j_2-j_3+m_1)} q^{-j_2(j_2+1) +j_3(j_3+1)
 -j_1(j_3+j_2+1)}q^{-{1\over 2} m_1}
 \cdot  \nonumber \\
&&\hskip -2 cm \cdot q^{-(j_2+j_3)m_1}
\sum_p { q^{2p(j_2+j_3 +1)} \over  [p]! [j_1+j_2-j_3-p]! [j_1+m_1-p]!
[j_3-j_2-m_1+p]! }
 \nonumber
 \end{eqnarray}

Comparing with eq.(\ref{Mj1}) we get eq.(\ref{M=3j}) with $j_2 = j(x)$ and $j_3
=j(x) +\sigma_1$ where $j(x)$ is
given by eq.(\ref{jx}). Moreover we find
\begin{eqnarray}
{\cal N}^{(j_1)}_\xi(m_1) &=& (-1)^{-{1\over 2}m_1 }q^{ {1\over 2}m_1 }
\nonumber \\
{\cal N}^{(j_1)}_\psi (x,\sigma_1) &=& (-1)^{-j_1+{1\over 2}(m_1-\sigma_1)}
{ \sqrt{ [j_1 +\sigma_1 ]![j_1 -\sigma_1 ]!}\over
 \prod_{r=1}^{j_1+\sigma_1} (1-x^2 q^{2r} ) }
{ (q-q^{-1})^{j_1}  x^{j_1} q^{j_1\sigma_1}
\over \Delta(j_1,j(x),j(x)+\sigma_1) \sqrt{[2j(x) +1 +\sigma_1 ]}}
\nonumber
\end{eqnarray}

\end{document}